\documentclass[10pt]{article}

\usepackage{dsfont}
\usepackage{graphics}
\usepackage{graphicx}
\usepackage{amssymb}
\usepackage{amsmath}
\usepackage{url}
\usepackage{color}

\newcommand{\rank}{\mathrm{Rank}}
\newcommand{\row}{\mathrm{Row}}
\newcommand{\col}{\mathrm{Col}}
\newcommand{\ODE}{\sc{ode}}

\def\proof{\textsc{Proof.} }
\def\foorp{\hfill$\square$}

 \newtheorem{thm}{Theorem}
 \newtheorem{lemma}{Lemma}
 
 \newtheorem{prop}{Proposition}
 \newtheorem{define}{Definition}
 \newtheorem{algorithm}{Algorithm}
 \newtheorem{Example}{EXAMPLE}[section]
 \newtheorem{remark}{Remark}
\begin{document}

\title{Exact Bivariate Polynomial Factorization in Q by
Approximation of Roots
\thanks{The work is partially supported by National Basic Research Program of China 2011CB302400
and the National Natural Science Foundation of
China(Grant NO.10771205).}}

\author{Yong Feng \ \ \ Wenyuan Wu\thanks{Corresponding author: dr.wenyuanwu@gmail.com}\ \ \ \
Jingzhong Zhang \vspace{0.2cm} \\
 \textit{Laboratory of Computer Reasoning and
Trustworthy Computing},
\\\textit{University of Electronic Sciences and Technology,}\\
\textit{Chengdu, P. R. China 611731} \\
}


\date{}
\maketitle

\begin{abstract}
Factorization of polynomials is one of the foundations of symbolic
computation. Its applications arise in numerous branches of
mathematics and other sciences. However, the present advanced
programming languages such as C++ and J++, do not support symbolic
computation directly. Hence, it leads to difficulties in applying
factorization in engineering fields. In this paper, we present an
algorithm which use numerical method to obtain exact factors of a
bivariate polynomial with rational coefficients. Our method can be
directly implemented in efficient programming language such C++
together with the GNU Multiple-Precision Library. In addition, the
numerical computation part often only requires double precision and
is easily parallelizable.

\end{abstract}

\section{Introduction}
Polynomial factorization plays a significant role in many problems
including the simplification, primary decomposition, factorized
Gr\"obner basis, solving polynomial equations and some engineering
applications, etc. It has been studied for a long time and some high
efficient algorithms have been proposed. There are two types of
factorization approaches. One is the traditional polynomial
factorization for exact input relying on symbolic computation, and
the other is approximate polynomial factorization for inexact input.

The traditional polynomial factorization methods  follow Zassenhaus'
approach \cite{Zassenhaus_a}\cite{Zassenhauss_b}. First,
Multivariate polynomial factorization is reduced to bivariate
factorization due to Bertini's theorem and Hensel lifting\cite{
GATHEN_a, GATHEN_b}. Then one of the two remaining variables is
specialized at random. The resulting univariate polynomial is
factored and its factors are lifted up to a high enough precision.
At last, the lifted factors are recombined to get the factors of the
original polynomial.

Approximate factorization is a natural extension of conventional
polynomial factorization. It adapts factorization problem to linear
algebra first, then applies numerical methods to obtain an
approximate factorization in complex which is the exact absolute
factorization of a nearby problem. In 1985, Kaltofen presented an
algorithm for computing the absolute irreducible factorization by
floating point arithmetic \cite{Kaltofen}. Historically, the concept
of approximate factorization appeared first in a paper on control
theory\cite{Mou}. The algorithm is as follows: 1) represent the two
factors $G$ and $H$ of the polynomials $F$ with unknown coefficients
by fixing their terms, 2) determine the numerical coefficients so as
to minimize $\|F-GH \|$. Huang et al \cite{Huang}. pursued this
approach, but it seems to be rarely successful, unless $G$ or $H$ is
a polynomial of several terms. In 1991, Sasaki et al. proposed a
modern algorithm\cite{Sasaki_a}, which use power-series roots to
find approximate factors. This algorithm is successful for
polynomials of small degrees. Subsequently, Sasaki et al. presented
another algorithm\cite{Sasaki_b} which utilizes zero-sum relations.
The zero-sum relations are quite effective for determining
approximate factors. However, computation based on zero-sum
relations is practically very time-consuming. In \cite{Sasaki_c},
Sasaki presented an effective method to get as many zero-sum
relations as possible by matrix operations so that approximate
factorization algorithm is improved. Meanwhile, Corless et al
proposed an algorithm for factoring bivariate approximate polynomial
based on the idea of decomposition of affine variety in
\cite{corless_B}. However, it is not so efficient to generalize the
algorithm to multivariate case. A major breakthrough is due to
Kaltofen et al. \cite{GKMYZ04,KMYZ08} who extended Gao's work
\cite{Gao03} from symbolics to numerics based on Ruppert matrix and
Singular Value Decomposition.

Symbolic factorization has been implemented in many Computer Algebra
System. However, it is difficult to implemented directly in
Programming Language such as C++ and J++, because most of
Programming Language standards do not support symbolic basic
operators, and the compilers do not implement the symbolic
computation, on which symbolic factorization is based. It restricts
exact factorization from being applied in many engineering fields.
Compared with symbolic factorization, approximate factorization can
be implemented more easily in the popular programming languages.
However it only gives approximate results even the input is exact.
In this paper, we propose an almost completely numerical algorithm,
which is not only implemented directly in the programming languages,
but also achieve exact results.

Except classic symbolic methods, some approaches have been proposed
to obtain exact output by approximation\cite{R.Kannan}\cite{zhang}.
The idea of obtaining exact polynomial factorization is from the
connect between an approximate root of a given polynomial and its
minimal polynomial in $Q$. Certainly, the minimal polynomial is a
factor of the given input. Based on lattice basis reduced algorithm
LLL and Integer Relation algorithm PSLQ of a vectors respectively,
there are two algorithms for finding exact minimal polynomial of an
algebraic number from its approximation. One is a numerical
algorithm \cite{R.Kannan,Qin} for factorization of a univariate
polynomial was provided by Transcendental Evaluation and high-degree
evaluation, and the other for factorization of bivariate polynomial
is based on LLL\cite{Hulst,Chen}. But they are not efficient.

In this paper, relying on LLL algorithm, we present an
almost-completely numerical method for exact factoring polynomial
with rational coefficient in $\mathbb{Q}$. First, we choose a sample
point in $\mathbb{Q}^{n-1}$ at random. After specialization  (i.e.
substitution) at the point, the roots of the resulting univariate
polynomial can be found very efficiently up to arbitrarily high
accuracy. Then applying minimal polynomial algorithm to these roots
yields an exact factorization of the univariate polynomial in
$\mathbb{Q}$. Next we shall move the sample point in ``good
direction" to generate enough number of points by using numerical
continuation. Especially, for the rest sample points, the
corresponding exact factorization can be found by using the same
combination of roots as found in the first step. And these roots
give more univariate polynomials for next step. Finally, the
multivariate factorization can be obtained by interpolation.

The paper is organized as follows. Section 2 gives a brief
introduction of the preparation knowledge. Minimal polynomial
algorithm will be discussed in Section 3. Then we present our method
in Section 4,5 and 6.

\section{Preparation}\label{sec:prep}

In this section, we briefly introduce the background knowledge and
related topics to our article.

\subsection{Homotopy Continuation Methods}

Homotopy continuation methods play a fundamental role in Numerical
Algebraic Geometry and provide an efficient and stable way to
compute all isolated roots of polynomial systems. These methods have
been implemented in many software packages e.g. Hom4PS
\cite{Hom4PS}, Bertini \cite{bertini}, PHCpack \cite{V99}.

The basic idea is to embed the target system into a family of
systems continuously depending on parameters. Then each point in the
parameter space corresponds to a set of solutions. Suppose we know
the solutions at a point. Then we can track the solutions from this
starting point to the point representing the target system we want
to solve.

First let us look at the simplest case: a univariate polynomial
$f(z)$ with degree $d$. We know that $f(z)$ has $d$ roots in
$\mathbb{C}$ (counting multiplicities). Of course we can embed
$f(z)$ into the family $a_{d} z^d+ a_{d-1}z^{d-1} + \cdots + a_0$,
where the $a_i$ are parameters. Now choose a start point
corresponding to $z^d -1 $ in this parameter space, whose roots are
\begin{equation}\label{eq:unity}
    z^0_k = e^{2k\pi \sqrt{-1}/d}, \; k= 0,1,...,d-1
\end{equation}
Then we use a real straight line in the parameter space to connect
$z^d -1$ with $f(z)$:
\begin{equation}\label{eq:uni-homotopy}
    H(z,t):= tf(z) + (1-t)(z^d-1).
\end{equation}
This form is a subclass of the family depending on only one real
parameter $t \in [0,1]$.

When $t=0$ we have the start system $H(z,0)=z^d-1$ and when $t=1$ we
have our target system $H(z,1)=f(z)$. An important question is to
show how to track individual solutions as $t$ changes from $0$ to
$1$. Let us look at the tracking of the solution $z_k$ (the $k$-th
root of $f(z)$). When $t$ changes from $0$ to $1$, it describes a
curve, which is function of $t$, denoted by $z_k= z_k(t)$. So
$H(z_k(t), t) \equiv 0$ for all $t \in [0,1]$. Consequently, we have
\begin{equation}\label{eq:uni-homotopy2}
    0 \equiv \frac{dH(z_k(t), t)}{dt} = \frac{\partial H(z,t)}{\partial
z} \frac{d z_k(t)}{dt} + \frac{\partial H(z,t)}{\partial t}.
\end{equation}
This problem is reduced to an {\ODE} for the unknown function
$z_k(t)$ together with an algebraic constraint $H(z_k(t), t) \equiv
0$. The initial condition is the start solution $z_k(0) = z^0_k$ and
$z_k(1)$ is a solution of our target problem $f(z)=0$.

\begin{remark}
In the book \cite{BCSS98}, Blum, Smale et al. show that on average
an approximate root of a generic polynomial system can be found in
polynomial time. Also application of the polynomial cost method for
numerically solving differential algebraic equations \cite{ICR05}
gives polynomial cost method for solving homotopies.
\end{remark}

But there is a prerequisite for the continuous tracking:
$\frac{\partial H(z,t)}{\partial z} \neq 0$ along the curve
$z=z_k(t)$. If the equations $z - z_k(t) =0$ and $tf'(z) + d (1-t)
z^{d-1} =0 $ have intersection at some point $(t,z_k(t))$, then we
cannot continue the tracking. There is way to avoid this singular
case, called the ``gamma trick" that was first introduced in
\cite{MS87}. We know two complex curves almost always have
intersections at complex points, but here $t$ must be real. So if we
introduce a random complex transformation to the second curve, the
intersection points will become complex points and such a
singularity will not appear when $t \in [0,1)$. Let us introduce a
random angle $\theta \in [-\pi,\pi]$ and modify the homotopy
(\ref{eq:uni-homotopy}) to
\begin{equation}\label{eq:uni-homotopy3}
    H(z,t):= t f(z) + e^{i\theta} (1-t)(z^d-1).
\end{equation}
It is easy to show that the $k$-th starting solution is still
$z^0_k$ in (\ref{eq:unity}) and that $z_k(1)$ is still a root of
$f(z)$.

\subsection{Genericity and Probability One} In an idealized model
where paths are tracked exactly and the random angle can be
generated to infinite precision, the homotopy
(\ref{eq:uni-homotopy3}) can be proved to succeed ``with probability
one". To clarify this statement, it is necessary to use a
fundamental concept in algebraic geometry: \textit{genericity}.

\begin{define}[Generic]
Let $X$ be an irreducible algebraic variety. We say a property $P$
holds \textbf{generically} on $X$, if the set of points of $X$ that
do not satisfy $P$ are contained in a proper subvariety $Y$ of $X$.
The points in $X\backslash Y$ are called \textbf{generic} points.
\end{define}

The set $X\backslash Y$ is called a Zariski open set of $X$. Roughly
speaking, if $Y$ is a proper subvariety of an irreducible variety
$X$ and $p$ is a random point on $X$ with uniform probability
distribution, then the probability that $p \notin Y$ is one. So we
can consider a random point as a generic point on $X$ without a
precise description of $Y$. Many of the desirable behaviors of
homotopy continuation methods rely on this fact.

\subsection{Coefficient-Parameter Homotopy}

There are several versions of the Coefficient-Parameter theorem in
\cite{SW05}. Here we only state the basic one.

\begin{thm}\label{thm:coeff-para}
Let $F(z;q)= \{f_1(z;q),...,f_n(z;q)\}$ be a polynomial system in
$n$ variables $z$ and $m$ parameters $q$. Let $\mathcal{N}(q)$
denote the number of nonsingular solutions as a function of $q$:
\begin{equation}\label{N(q)}
    \mathcal{N}(q):= \# \left\{z \in \mathbb{C}^n : F(z;q) = 0, \; \det \left(\frac{\partial F}{\partial z}(z;q) \right) \neq 0 \right\}
\end{equation}
Then,
\begin{enumerate}
    \item There exist $N$, such that $\mathcal{N}(q) \leq
    N$ for any $q \in \mathbb{C}^m$. Also $\{q \in \mathbb{C}^m: \mathcal{N}(q) = N \}$ is a
    Zariski open set of $\mathbb{C}^m$. The exceptional set $Y = \{q: \mathcal{N}(q) <
    N\}$ is an affine variety contained in a variety
    with dimension $m-1$.
    \item The homotopy $F(z;\phi(t)) = 0$ with $\phi(t) : [0,1) \rightarrow \mathbb{C}^m \backslash
    Y$ has $N$ continuous non-singular solution paths $z(t)$.
    \item When $t \rightarrow 1^-$, the limit of $z_k(t), \; k=1,. .., N$ includes all
    the non-singular roots of $F(z;\phi(1))$.
\end{enumerate}
\end{thm}

An important question is how to choose a homotopy path $\phi(t)$
which can avoid the exceptional set $Y$. The following lemma
\cite{SW05} gives an easy way to address this problem.

\begin{lemma}\label{lem:line-homo}
Fix a point $q$ and a proper algebraic set $Y$ in $\mathbb{C}^m$.
For a generic point $p \in \mathbb{C}^m$, the one-real-dimensional
open line segment $\phi(t):= (1-t) \;p + t \; q, t \in [0,1)$ is
contained in $\mathbb{C}^m \backslash Y$.
\end{lemma}

\subsection{Reductions}
Before factorization of a given polynomial, we shall first apply
certain reductions to the input to obtain a square-free polynomial
over $\mathbb{Q}$, which can remove multiplicities and ease the
computation of the roots. Also we can assume each factor involves
all the variables and has more than one term. Otherwise, we can
compute the GCD to reduce the problem. For example, let $F= f(x,y)
g(y)$. Then $F_x = f_x g$ and $\gcd(F,F_x) = g$ which gives us the
factor $g(y)$.

By the Hilbert Irreducibility Theorem, we can further reduce the
problem to univariate case by generic (random) specialization of one
variable to a rational number. More precisely, if $f(x,y)$ is
irreducible in $\mathbb{Q}[x,y]$, then for a generic rational number
$y_0$, $f(x,y_0)$ is also irreducible in the ring $\mathbb{Q}[x]$.
It means that the factorization is commutable with generic
specialization.

For univariate polynomial, there are symbolic methods to preform
exact factorization in $\mathbb{Q}$. Here we are more interested in
numerical methods, namely from approximate roots to exact factors.

\section{Minimal Polynomial by Approximation}\label{sec:Minipoly}

There are two methods to compute the minimal polynomial of an
algebraic number from its approximation. One is based on the LLL
algorithm of the basis reduction\cite{R.Kannan}, and another is
based on PSLQ\cite{Qin}. The later one is more efficient than the
former one. However, it can only compute the minimal polynomial of a
real algebraic number while the former one can find minimal
polynomial of a complex algebraic number. Hence, we introduce the
former algorithm which is more suitable for this paper here. We
refer the reader to the paper \cite{R.Kannan} for more details.

Let $p(x)=\sum_{i=0}^{i=n}p_ix^i$ be a polynomial. The length $|p|$
of $p(x)$ is defined as the Euclidean norm of the vector
$(p_0,p_1,\cdots,p_n)$, and the height $|p|_{\infty}$ as the
$L_{\infty}$-norm of the vector $(p_0,p_1,\cdots,p_n)$. The degree
and height of an algebraic number are defined as the degree and
height, respectively, of its minimal polynomial.

Suppose that we have upper bound $d$ and $H$ on the degree and
height respectively of an algebraic number with $|\alpha|\le 1 $,
and a complex rational number $\bar{\alpha}$ approximating $\alpha$
such that $|\bar\alpha|\le 1$ and $|\alpha-\bar\alpha|<2^{-s}/(4d)$,
where $s$ is the smallest positive integer such that
$$2^s>2^{d^2/2}(d+1)^{(3d+4)/2}H^{2d} $$
\newcounter{num}

\begin{algorithm} \label{alg:miniPoly}[miniPoly]

For $n=1,2,\cdots,d$ in succession, do the following steps
\begin{list}{Step \arabic{num}:}{\usecounter{num}\setlength{\rightmargin}{\leftmargin}}
\item construct
\begin{equation}\label{ls_matrix}
\left(
\begin{array}{ccccccc}
1&0&0&\cdots&0&2^s\cdot Re(\bar\alpha_0)&2^s\cdot Im(\bar\alpha_0)\\
0&1&0&\cdots&0&2^s\cdot Re(\bar\alpha_1)&2^s\cdot Im(\bar\alpha_1)\\
0&0&1&\cdots&0&2^s\cdot Re(\bar\alpha_2)&2^s\cdot Im(\bar\alpha_2)\\
\vdots&\vdots&\vdots&\ddots&\vdots&\vdots&\vdots\\
0&0&0&\cdots&1&2^s\cdot Re(\bar\alpha_n)&2^s\cdot Im(\bar\alpha_n)
\end{array}\right)
\end{equation}
where $Re(a)$ and $Im(a)$ stand for the real part and imaginary
part, respectively, of complex $a$, $\alpha_0=1$ and
$|\bar{\alpha_i}-\bar{\alpha}^i|\le 2^{-s-1/2}$ for
$i=1,2,\cdots,n$. Note $\bar{\alpha_i}$ can be computed by rounding
the powers of $\bar\alpha$ to $s$ bits after the binary points.
\item Denote by $b_i$ the row $i+1$ of the matrix in
(\ref{ls_matrix}). Apply the basic reduction algorithm to lattice
$L_s=(b_0,b_1,\cdots,b_n)$, and obtain the reduced basis of the
lattice.
\item If the first basis vector $\tilde{v}=(v_0,v_1,\cdots,v_n,v_{n+1},v_{n+2})$ in the reduced basis
satisfies $|\tilde{v}|\le 2^{d/2}(d+1)H$, then return polynomial
$v(x)=\sum_{i=0}^nv_ix^i$ as the minimal polynomial of algebraic
number $\alpha$.
\end{list}
\end{algorithm}

Note: It is no major restriction to consider $\alpha$ with
$|\alpha|\le 1$ only. In fact, if $|\alpha|>1$ satisfies the
polynomial $h(x)=\sum_{i=0}^dh_ix^i$, then $1/\alpha$ satisfies the
polynomial $\sum_{i=0}^dh_{d-i}x^i$. Therefore, if
$\sum_{i=0}^dh_{d-i}x^i$ is computed, the $h(x)$ is obtained.
Furthermore, an $\varepsilon$-approximation $\bar\alpha$ to $\alpha$
with $|\alpha|>1$ easily yields a $3\varepsilon$-approximation
$\bar\beta$ to $\beta=1/\alpha$. This can be easily verified.

The following theorem shows the computation amount of calculating
the minimal polynomial of an algebraic number\cite{R.Kannan}:
\begin{thm}\label{thm:minipoly}
Let $\alpha$ be an algebraic number, and let $d$ and $H$ be upper
bounds on the degree and height, respectively, of $\alpha$. Suppose
that we are given an approximation $\bar\alpha$ to $\alpha$ such
that $|\alpha-\bar\alpha|\le 2^{-s}/(12d)$, where $s$ is the
smallest positive integer such that
$$2^s>2^{d^2/2}(d+1)^{(3d+4)/2}H^{2d} $$
Then the minimal polynomial of $\alpha$ can be determined in
$O(n_0\cdot d^4(d+\log H))$ arithmetic operations on integers having
$0(d^2\cdot(d+\log H))$ binary bits, where $n_0$ is the degree of
$\alpha$.
\end{thm}

\section{Finding More Polynomials by Continuation}\label{sec:FM}

In the previous stage, we have the factors after specialization ,
which are univariate polynomials. To construct the factor of two
variables by using interpolation, we need more information, i.e.
specializations at more points. The main tool is the homotopy
continuation method.

\subsection{Applying numerical continuation to factorization}
Suppose an input polynomial $F(x,y)$ is reducible. Geometrically, if
$\mathcal{C}$ denotes the zero set of $f$ i.e. the union of many
curves in $\mathbb{C}^2$, removing the singular locus of
$\mathcal{C}$ from each curve $\mathcal{C}_i$, the regular sets
$\mathcal{S}_i$ are connected in $\mathbb{C}^2$. Moreover, the
singular locus has lower dimension, consequently it is a set of
isolated points.

Suppose $f(x,y)$ is an irreducible factor of $F$ in $\mathbb{Q}$.
Let $y_0,y_1$ be random rational numbers.
 By the Hilbert Irreducibility Theorem the univariate polynomials $f_0=f(x,y_0)$ and
 $f_1=f(x,y_1)$ are irreducible as well. Suppose we know the roots
 of  $f_0$. Then we can choose a path to connect $y_0$ and $y_1$ avoiding the singular
 locus which has measure zero. By the Coefficient-Parameter Theorem,
 all the roots of $f_1$ can be obtained by the following homotopy continuations:
\begin{equation}\label{eq:homotopy1}
      \left\{ \begin{array}{cl}
              f(x,y)=0  &  \\
              (1-t)(y-y_0) + t(y-y_1)\gamma =0 & \\
              \end{array}\right.
\end{equation}
Moreover any generic complex number $\gamma$ implies that the
homotopy path can avoid the singular locus by Lemma
\ref{lem:line-homo} when we track the path.

\subsection{Control of the precision}
Let $\{x_1,..,x_m\}$ be the exact roots of $f_1$ and $g$ be the
primitive polynomial of $f_1$. Then
\begin{equation}\label{eq:prod_roots}
g= \alpha \prod_{i=1}^m (x-x_i) \in \mathbb{Z}[x],
\end{equation} for some
integer number $\alpha$.

Note that we only have the approximate roots
$\{\tilde{x}_1,..,\tilde{x}_m\}$.

\begin{prop}\label{prop:root2poly}
Let $p= \prod_{i=1}^m (x-x_i)$ and  $\tilde{p}= \prod_{i=1}^m
(x-\tilde{x}_i)$. Let $\delta = \max_{i=1,..,m}\{|x_i-
\tilde{x}_i|\}$ and $r= \max_{i=1,..,m}\{|\tilde{x}_i|\}$. If
$\delta$ is sufficiently small. Then
\begin{equation}
||\tilde{\textbf{p}}- \textbf{p}||_{\infty}\leq
\left(\max_{i=1,..,m}\{ \;r^{i-1} {{m-1}\choose{i-1}} \}\;m +1
\right) \; \delta
\end{equation}
\end{prop}
\proof Let $x_i = \tilde{x}_i + \delta_i$. Thus, $| \delta_i| \leq
\delta$. The left hand side $||\tilde{\textbf{p}}-
\textbf{p}||_{\infty} = ||\prod_{i=1}^m (x-x_i+\delta_i) -
\prod_{i=1}^m (x-x_i)|| = || \sum_{j=1}^m \prod_{i\neq
j}(x-x_j)\delta_j || + o(\delta)$. An upper bound of the
coefficients of $\prod_{i\neq j}(x-x_j)$ with respect to $x^{m-i}$
is ${{m-1}\choose{i-1}}r^{i-1}$. Hence, $||\tilde{\textbf{p}}-
\textbf{p}||_{\infty}\leq  \max_{i=1,...,m}
{{m-1}\choose{i-1}}r^{i-1} \cdot m \delta + \delta$ \foorp

Now let us consider how to find $\alpha$. Suppose the input
polynomial is $F(x,y)$ and $f$ is a factor of $F$. The primitive
polynomial of $f(x,y_1)$, which is $g$, must be a factor of the
primitive polynomial of $F(x,y_1)$. Thus, the leading coefficient of
$g$ must be a factor of the leading coefficient of the primitive
polynomial of $F(x,y_1)$. Let $\alpha$ be the leading coefficient of
the primitive polynomial of $F(x,y_1)$. Then let $p= \alpha
\prod_{i=1}^m (x-x_i) \in \mathbb{Z}[x]$. Note that itself may not
be primitive, but its primitive polynomial is $g$.

Let $M =\max_{i=1,..,m}\{i \;r^{i-1} {{m}\choose{i}} \} +1 $.
$\tilde{p} = \alpha \prod_{i=1}^m (x-\tilde{x}_i)$. Thus, if $\delta
< \frac{1}{2\alpha M}$ then $||p-\tilde{p}||_{\infty} < 0.5$. It means that
we can round each coefficient of $\tilde{p}$ to the nearest integer
to obtain exact polynomial $p$ which gives $g$.

\subsection{Detecting the degrees of factors}

After specialization at $y=y_0$, we obtain the information about the
number of factors and the degree of each factor with respect to $x$.
The degrees with respect to $y$ of factors provide the bound of the
number of interpolation nodes. Certainly, we can use the degree of
the input $\deg_{y}(F)$ as the bound. However, the degrees  with
respect to $y$ of factors are usually much less than $\deg_{y}(F)$,
especially when there are many factors. Therefore, for high
efficiency, it is better to know the degree with respect to $y$ of
each factor. Now we will apply numerical continuation to detect such
degree information.

We define the notation of $2$-tuple degree to be $$ \deg(f) =
[\deg_{x}(f),  \deg_{y}(f)].$$

Suppose $\deg(f) = [m,n]$ and $f$ has $r$ factors. Applying an
approach of univariate polynomial solving to $f(x,y_0)$ and
$f(x_0,y)$ yields points on the curve $A = \{(x_1, y_0), (x_2,
y_0),...,(x_m, y_0) \}$ and $B = \{(x_0, y_1), (x_0, y_2),...,(x_0,
y_n) \}$ respectively. In addition, we also know the decomposition
of two points sets in $r$ groups with cardinalities
$\{a_1,...,a_r\}$ and $\{b_1,...,b_r\}$. Moreover $\sum a_i = m$ and
$\sum b_i = n$.

Choose one point from each group of the first set $A$. Starting from
these points, we track the homotopy path
\begin{equation}\label{eq:homotopy2}
      \left\{ \begin{array}{cl}
              f(x,y)=0  &  \\
              (1-t)(y-y_0) + t(x-x_0)\gamma =0 & \\
              \end{array}\right.
\end{equation}
Because of the genericity of the choice of $y_0$, $x_0$ and
$\gamma$, the path avoids the singular locus. When $t=1$, the
endpoint must belong to the second set $B$. For example if the
starting point of the first group of $A$ and its end point belongs
to the $i$th group of $B$. Then we know the first factor has degree
$[a_1, b_i]$. Similarly, the degrees of other factors can be
detected in the same way.

\section{Interpolation}\label{sec:Interp}

Polynomial interpolation is a classical numerical method. It is
studied very well for univariate polynomials in numerical
computation. Polynomial interpolation problem is to determine a
polynomial $f(x)\in F[x]$ with degree not greater than $n\in
\mathbb{N}$ for a given pairs $\{(x_i,f_i),i=0,\cdots,n\}$
satisfying $f(x_i)=f_i$ for $i=0,\cdots,n$, where $F$ is a field and
$x_i,f_i\in F$. In general, there are four types of polynomial
interpolation method: Lagrange Interpolation, Neville's
Interpolation, Newton's Interpolation and Hermite Interpolation.
Lagrange interpolation and  Newton's Interpolation formula are
suited for obtaining interpolation polynomial for a given set
$\{(x_i,f_i), i=0,\cdots,n\}$.  Neville's interpolation method aims
at determining the value of the interpolating polynomial at some
point.  If the interpolating problem prescribes at each
interpolation point $\{x_i,i=0,\cdots,n\}$ not only the value but
also the derivatives of desired polynomial, then the Hermite formula
is preferred.

Different from the traditional interpolation problem above, our
problem is to construct a bivariate polynomial from a sequence of
univariate polynomials at chosen nodes. It is important to point out
that the univariate polynomials are constructed by roots, which may
not be equal to the polynomials by substitutions. But the only
difference for each polynomial is just a scaling constant.

More precisely, in this paper, we aim to solve a special polynomial
interpolation problem: \textit{given a set of nodes and square free
polynomials $\{(y_i\in F,f_i(x)\in F[x]),i=0,\cdots,k\}$, compute a
square free polynomial $f(x,y)\in F[x,y]$ of degree with respect to
$x$ not greater than $n$, where $F$ is a field, such that $f(x,y_i)$
and $f_i(x)$ have the same roots.}


\subsection{Illustrative examples}

\begin{Example}
Let $f = x^2+y^2 -1$. Since its degree with respect to $y$ is two,
we need three interpolation nodes which are $y=-1/2, 0, 1/2$.
Suppose we know the roots at each node, then the interpolating
polynomials are $\{f_0= x^2 -3/4, f_1= x^2 -1, f_2= x^2 -3/4 \}$. To
construct original polynomial $f$, we can use Lagrange method.

Let $\ell_1 = \frac{y(y-1/2)}{(-1/2-0)(-1/2-1/2)} = 2y^2-y$.
Similarly, $\ell_2 = -4y^2+1$ and $\ell_3 = 2y^2+y$. It is easy to
check that $(x^2 -3/4)\ell_1 + (x^2 -1)\ell_2 + (x^2 -3/4)\ell_3 =
f$.
\end{Example}

In the example above, the coefficient of $f$ with respect to $x^2$
is a constant $1$. Making the interpolating polynomials given by
(\ref{eq:prod_roots}) monic, we can construct $f$ correctly by
Lagrange basis. However, if the coefficient is nonconstant, i.e. a
polynomial of $y$, then it is not straightforward to find $f$. The
example below shows this problem.

\begin{Example} \label{example:2}
Let $f = xy -1$. The nodes are $y=2,3$. We know the roots are $1/2,
1/3$ respectively at the nodes. Then the monic interpolating
polynomials are $\{x-1/2, x-1/3\}$.  If we still apply Lagrange
basis $\ell_1 = -y+3, \; \ell_2 = y-2$, it gives $(x-1/2)(-y+3) +
(x-1/3)(y-2)= x + 1/6 \;y -5/6$ which is totally different from the
target polynomial $xy -1$.
\end{Example}

The basic reason is that the interpolating polynomials are not the
polynomials after specialization s, and the only difference is
certain scaling constants. To find these constants, we need more
information. Now we use one more node: when $y=4$, the monic
interpolating polynomial is $x-1/4$. By multiplying a scaling
constant to $f$ we can assume $f(x,4) = x-1/4$, then there exist
$a,b$ such that $f(x,2) = a(x-1/2)$ and $f(x,3) = b(x-1/3)$. The
corresponding Lagrange bases are $\ell_1 = (y-3)(y-4)/2, \; \ell_2 =
-(y-2)(y-4),\;  \ell_3 = (y-2)(y-3)/2$. Then $f = a(x-1/2)\ell_1 +
b(x-1/3)\ell_2 + (x-1/4)\ell_3$. The coefficient of $f$ with respect
to $y^2$ must be zero. Consequently we have
\begin{equation}
1/2\, \left( x-1/2 \right)a + \left( x-1/3 \right)b +1/2\,x-1/8 = 0
\end{equation}
which implies a linear system $$ 1/2 \; a-b+1/2 = 0,
-1/4\;a+1/3\;b-1/8 = 0$$ The solution is $a=1/2, b=3/4$.
Substituting  them back to two nodes interpolation formula yields
the polynomial we need, up to a constant $1/4$
$$1/2(x-1/2)(-y+3) + 3/4(x-1/3)(y-2)= (xy -1)/4$$

\subsection{Interpolation with indeterminates}

To extend the idea in example \ref{example:2}, we present a method to construct
desired bivariate polynomial by using monic univariate interpolating
polynomials.

Suppose $f$ is irreducible and its degrees with respect to $x$
and $y$ are $m$ and $n$ respectively. Consider $x$ as the main
variable, we can express this polynomial by $f= \sum_{i=0}^m c_i(y)
x^i$, where $c_i$ are polynomials of $y$ of degree less than or
equal to $n$. We can consider each $c_i$ as a vector in monomial
basis. Suppose there are $r$ linearly independent coefficients. If
$r=1$, then $c_i(y) = a_i c_0(y)$ for some constant $a_i$ and $f =
(\sum_{i=0}^m a_i x^i) \cdot c_0(y)$. It contradicts the assumption
that $f$ is irreducible. Hence, $r\geq 2$.

Now we consider how to construct $f$ by using the interpolating
polynomials $\{f_0(x),f_1(x),...,f_k(x)\}$ at $k+1$ nodes $\{y_0,
y_1,..., y_k\}$ respectively chosen at random.

Let $C$ be a $(k+1)\times (m+1)$ matrix  $[\textbf{c}_0, ...,
\textbf{c}_m]$ where $\textbf{c}_i$ is the column vector in monomial
basis $\{y^k,y^{k-1},...,1\}$ of the polynomial $c_i$. Let $V$ be
the Vandermonde matrix $\left(
           \begin{array}{cccc}
             y_0^k & y_0^{k-1} & \cdots & 1 \\
             y_1^k & y_1^{k-1} & \cdots & 1 \\
             \vdots & \vdots & \vdots & \vdots \\
             y_k^k & y_k^{k-1} & \cdots & 1 \\
           \end{array}
         \right)$.
Let $A$ be a $(k+1) \times (m+1)$ matrix where $A_{ij}$ is the
coefficient of the $i$th interpolating polynomial with respect to
$x^j$. To make the solution unique, we may fix $f(x,y_0) = f_0$ and
suppose $f(x,y_i) = \lambda_i f_i$ and $\lambda_i \neq 0$  for
$i=1,..,k$. Let $\Lambda = \left(
                   \begin{array}{cccc}
                     1 &  &  &  \\
                      & \lambda_1 &  &  \\
                      &  & \ddots &  \\
                      &  &  & \lambda_k \\
                   \end{array}
                 \right)$.

Therefore,
\begin{equation}\label{eq:lambda}
 V \cdot C = \Lambda \cdot A.
\end{equation}

Since $\{y_i\}$ are distinct, the Vandermonde matrix has inverse and
consequently $C = V^{-1} \cdot \Lambda \cdot A$. By our assumption,
the degree with respect to $y$ is $n$. It means that the first $k-n$
rows of $C$ must be zero. The zero at the $i$th row and $j$th column
corresponds an equation. Thus, it leads to a linear system
\begin{equation}\label{eq:scalingConst}
\row(V^{-1}, i) \cdot \Lambda \cdot \col(A,j) = 0,
\end{equation}
for $1\leq i \leq k-n$ and  $1 \leq j  \leq m+1$ with $k$ unknowns.

Only $r$ linearly independent columns in $A$, so there are $(k-n) \;
r$ equations and $k$ unknowns. The existence of the solution is due
to the origination of the interpolating polynomials $f(x,y_i) =
\lambda_i f_i$ for $i=1,...,k$. The linear system has unique
solution implies that $(k-n)r \geq k$. Thus, $k \geq rn/(r-1)$.  Let
$\mu$ be the smallest integer greater than or equal to
$\frac{rn}{r-1}$, namely
\begin{equation}\label{eq:k}
\mu = \lceil \frac{rn}{r-1} \rceil.
\end{equation}

Thus, to determine the scaling constants $\{ \lambda_i \}$, we need
at least $\mu$ more interpolation nodes.

To find an upper bound for the number of nodes, let us consider $f$
as a monic polynomial with rational function coefficients. All the
coefficients $\{c_m,...,c_0\}$ can be uniquely determined by
rational function interpolation of $ x^m + c_{m-1}/c_m x^{m-1} +
\cdots + c_0/c_m$ at $2n+1$ nodes. Therefore, it requires $2n$ nodes
except the initial one. Thus, we have $\mu \leq k \leq 2n$.

But this upper bound is often overestimated, and for some special
case the polynomial $f$ can be constructed by using less nodes.

\begin{prop}
Let $f$ be a polynomial in $\mathbb{Q}[x,y]$ and  $\deg(f)=[m,n]$.
Suppose $m \geq n$ and f has $n+1$ linearly independent
coefficients.  Then $f$ can be uniquely determined by $n+2$ monic
interpolating polynomials $\{f_0(x),...,f_{n+1}(x)\}$ up to a
scaling constant.
\end{prop}
\proof Suppose the first $n+1$ columns of $A$ are linearly
independent. By Equation (\ref{eq:lambda}), we construct $n+1$
equations: $ \row(V^{-1}, 1) \cdot \Lambda \cdot \col(A,j) = 0$, for
$j=1,...,n+1$. Let $B$ be the transpose of the submatrix consisting
of the first $n+1$ columns of $A$ and $\textbf{v} =
(v_1,...,v_{n+2})^t$ be the transpose of $\row(V^{-1}, 1)$. Thus, \\
$$\textbf{0} = B \cdot \left(
                   \begin{array}{cccc}
                     \lambda_1 &  &  &  \\
                      & \lambda_2 &  &  \\
                      &  & \ddots &  \\
                      &  &  & \lambda_{n+2} \\
                   \end{array}
                 \right) \cdot \textbf{v} = B \cdot
\left(
                   \begin{array}{cccc}
                     v_1 &  &  &  \\
                      & v_2 &  &  \\
                      &  & \ddots &  \\
                      &  &  & v_{n+2} \\
                   \end{array}
                 \right) \cdot (\lambda_1,...,\lambda_{n+2})^t$$
Because the first $n+1$ coefficients of $f$ are linearly
independent, the evaluations of them at $n+2$ random points must be
linearly independent. So the rank of $B$ is $n+1$. Here $v$ can be
expressed by explicit form of the Vandermonde matrix \cite{Turner66}
which is a vector of polynomials of $\{y_0,...,y_{n+1}\}$. For
generic choice of $\{y_0,...,y_{n+1}\}$, each $v_i \neq 0$. Hence,
the rank of $B \cdot \left(
                   \begin{array}{cccc}
                     v_1 &  &  &  \\
                      & v_2 &  &  \\
                      &  & \ddots &  \\
                      &  &  & v_{n+2} \\
                   \end{array}
                 \right)$ is still $n+1$ and its nullity equals one.
We can choose any solution $\{\lambda_i\}$ to construct $f$ by
Lagrange basis: $\sum_{i=0}^n  \lambda_i f_i \ell_i$. \foorp

\begin{remark}
In our algorithm, we compute $\{ \lambda_i \}$ starting from $\mu$
more nodes (together with the initial node $y_0$), and we add
incrementally more nodes if necessary. But interestingly, the
experimental results show that the lower bound $\mu$ is often
enough. This fact deserves further study.
\end{remark}

\begin{algorithm}\label{alg:NIF}[Interpolation]
\begin{tabbing}
\ \ \ \= Input :  a set of polynomials $\{f_0(x),...,f_k(x)\}
\subset \mathbb{Z}[x,y]$\\
\ \ \hspace{1cm} a set of rational numbers $\{y_0,...,y_k\}$ \\
\ \ \hspace{1cm} an integer $n$ the degree of $f$ with respect to $y$ \\
\> Output: $f \in \mathbb{Z}[x,y]$, such that $f(x,y_i) =  f_i(x)$.
\end{tabbing}
\vspace{-14pt}

\begin{enumerate}
  \item Let $A$ be the matrix consisting of the coefficient row vectors
  of the input univariate polynomials.
  \item Let $r = \rank(A)$. If $k < \mu$,
  then it needs more interpolation nodes.
  \item Solve the homogenous linear system (\ref{eq:scalingConst}) to obtain the scaling constants $\{\lambda_1,..,\lambda_k \}$
  \item If the solution is not unique, then it needs more interpolation nodes.
  \item Else $f = \sum_{i=0}^k  \lambda_i f_i \ell_i  \in \mathbb{Q}[x,y]$
  \item Return the primitive polynomial of $f$
\end{enumerate}

\end{algorithm}

\section{Combination of Tools}\label{sec:Algo}

Now we combine the tools introduced in previous sections to obtain a
new factorization algorithm. A preliminary version of the algorithm
is implemented in Maple. For the efficiency, it requires a more
sophisticate version in C++, even parallel program.

\subsection{Main steps of the algorithm}

\begin{algorithm}\label{alg:NIF}[Factorization]

$F$ =  \textbf{Fac}$(f)$
\vspace{-12pt}
{\rm
\begin{tabbing}
\ \ \ \= Input : \= $f$,  a primitive polynomial $f \in
\mathbb{Z}[x,y]$ such that $\gcd(f,f_x)=1$\\
\> Output: \= $F$,  a set of primitive polynomials $\{f_1,...,f_r\}
\subset \mathbb{Z}[x,y]$, such that $f = \prod f_i$.
\end{tabbing}
\vspace{-14pt}

\begin{enumerate}
  \item Apply a numerical solver to approximate the roots of $f(x,y_0)=0$ and $f(x_0,y)=0$ at generic points $x_0,y_0 \in \mathbb{Q}$.
  \item Apply \textbf{miniPoly} to roots above and decompose the solutions.
  And generate the minimal polynomials for them and we have grouping information of roots. In this step, it needs
  Newton iteration to refine the roots up to desired accuracy.
  \item Apply homotopy (\ref{eq:homotopy2}) to obtain the degrees of each factor.
  \item For group $i$ (corresponding to the factor $f_i$), $i=1,...,r$, use homotopy (\ref{eq:homotopy1}) to generate its approximate roots
  at random rational numbers $\{y_1,...,y_k\}$.
  \item For each set of roots at $y_j$, refine the roots to the accuracy given by Proposition \ref{prop:root2poly}, then make the product and
  construct the polynomial $f_i(x,y_j)$.
  \item Call \textbf{interpolate} with the interpolating polynomials
  $\{f_i(x,y_0),...,f_i(x,y_k)\}$ to construct $f_i(x,y)$.
\end{enumerate}
}

\end{algorithm}

\subsection{A simple example}
Let us consider a polynomial $f = (x\,y -2)\,(x^2+y^2-1)$. First, we
choose a sequence of random rational numbers $\{97/101, 1, 104/101,
123/101, 129/101,...\}$. Substituting $y=97/101$ into $f$ yields
Mignotte bound of the coefficients of factors $9170981$ and the
digits required to produce the minimal polynomial is $110$ by
Theorem \ref{thm:minipoly}.  Then compute the approximate roots of
$f(x,97/101)$ up to $110$ digits accuracy. The miniPoly subroutine
gives two groups of points: $[[1,2],[3]]$ and the corresponding
minimal polynomials $[-792+10201\,x^2, -202+97\,x]$. By Hilbert
Irreducibility theorem, there are two factors. On the other hand,
fix the value of $x$ and obtain the univariate polynomials
$[-202+97\,y, -792+10201\,y^2 ]$ and $[[3],[1, 2]]$.

Starting from the first point of group one, the Homotopy
(\ref{eq:homotopy2}) path ends at a point which satisfies $y,
-792+10201\,y^2$. It implies that $-792+10201\,x^2$ and $y,
-792+10201\,y^2$ are from the same factor of degree $[2,2]$. By
Equation (\ref{eq:k}), we need $\mu = \lceil \frac{rn}{r-1}\rceil =
4$ more interpolating polynomials which are produced by Homotopy
(\ref{eq:homotopy1}). Thus, there are five polynomials
$[-792+10201\,x^2, x^2, 615+10201\,x^2, 4928+10201\,x^2,
6440+10201\,x^2]$. The scaling constants $[\lambda_1 = 1, \lambda_2
= 10201, \lambda_3 = 1, \lambda_4 = 1, \lambda_5 =1]$ are obtained
by system (\ref{eq:scalingConst}). Consequently, the Lagrange
interpolation formula gives the correct factor $-1+x^2+y^2$.

Since the degree of the other factor is $[1,1]$, it needs $\mu = 2$
more polynomials and they are $[-202+97\,x, x-2, -101+52\,x]$. The
corresponding scaling constants are $[\lambda_1 = 1/2, \lambda_2 =
101/2, \lambda_3 = 1]$ and the resulting factor is $x\,y-2$.

\section{Conclusion}

A new numerical method to factorize bivariate polynomials exactly is
presented in this article. We implemented our algorithm in Maple to
verify the correctness. More importantly, the main components of our
algorithm, miniPoly and numerical homotopy continuation can be
implemented directly in C++ or J++ with existing multi-precision
packages, e.g. GNU MP library. Furthermore, these two numerical
components are naturally parallelizable. Therefore, it gives an
alternative way to exact factorization which can take the advantages
of standard programming languages and parallel computation
techniques widely used by industries.

In this article, we mainly focus on bivariate case. It is quite
straightforward to extend to multivariate case. However, the number
of the interpolation nodes grows exponentially as the increasing of
the number of monomials. A more practical way to deal with such
difficulty is to exploit the sparsity if the factors are sparse. It
desires the further study in our future work.

\end{document}